\documentclass[11pt]{amsart}
\usepackage{amsmath,amscd,amssymb,amsfonts}
\theoremstyle{plain}
\newtheorem{thm}{Theorem}
\newtheorem{lem}[thm]{Lemma}
\newtheorem{cor}[thm]{Corollary}

\numberwithin{thm}{section}
\numberwithin{equation}{section}

\newcommand{\eq}[2]{\begin{equation}\label{#1}#2 \end{equation}}

\newcommand{\gr}{{\rm gr}}

\newcommand{\Pic}{{\rm Pic}}

\renewcommand{\div}{{\rm div}}

\newcommand{\sD}{{\mathcal D}}

\newcommand{\sO}{{\mathcal O}}

\newcommand{\sV}{{\mathcal V}}

\newcommand{\C}{{\mathbb C}}

\renewcommand{\P}{{\mathbb P}}
\newcommand{\Q}{{\mathbb Q}}

\newcommand{\Z}{{\mathbb Z}}

\begin{document}

\title[The transcendental $K_1$ of a $K3$ surface]{The transcendental part of 
the regulator map for $K_1$ on a mirror family of $K3$ surfaces}
 
\author{Pedro Luis del Angel}
\address{Pedro Luis del Angel, Cimat, Guanajuato, Mexico}
\email{luis@cimat.mx}

\author{Stefan M\"uller-Stach}
\address{Stefan M\"uller-Stach, Fachbereich Mathematik, 45117 Essen, Germany}
\email{mueller-stach@uni-essen.de}

\date{}
\begin{abstract} We compute the transcendental part of the 
normal function corresponding to the Deligne class of a cycle in 
$K_1$ of a mirror family of quartic $K3$ surfaces. 
The resulting multivalued function does not satisfy
the hypergeometric differential equation of the periods and we conclude that
the cycle is indecomposable for most points in the mirror family. The
occuring inhomogenous Picard-Fuchs equations are related to
Painlev\'e VI type differential equations. 
\end{abstract}
\thanks{We are grateful to DFG, UAM and Univ. Essen for supporting this project}
\maketitle

\section{The regulator map and Picard-Fuchs equations}

In this paper we study the first non-classical higher K-group
$K_1(X)$ for a smooth complex projective surface $X$. It was conjectured
by H. Esnault around 1995 that certain elements in this group can be detected
in the transcendental part of the Deligne cohomology group $H^3_{\sD}(X,\Z(2))$
via the regulator (Chern class) map. The transcendental part of the regulator map is
defined as an Abel-Jacobi type integral of holomorphic two-forms over non-closed real
2-dimensional chains in $X$ associated to these elements. At that time is was only
known that one could detect such classes in the complementary $(1,1)$-part of Deligne
cohomology (see e.g. \cite{SMS}). The goal of our paper is to show that Esnault's
conjecture is true by looking at the differential equations which are satisfied by
the normal functions arising from such classes in a family of surfaces. It turns out
that the resulting equations for Abel-Jacobi type integrals with parameters are
strongly connected to a generalization of Painlev\'e VI type differential equations.
\\
The higher K-groups $K_1(X),K_2(X),\ldots $ of an algebraic variety $X$ 
were defined around 1970 by D. Quillen \cite{Q}. 
Later Bloch \cite{B2} showed that on smooth quasi-projective varieties 
all their graded pieces with respect to the 
$\gamma$-filtration may be computed as   
$$
\gr_\gamma^p K_n(X)_\Q \cong CH^p(X,n)_\Q
$$
where $CH^p(X,n)$ are Bloch's higher Chow groups \cite{B2}. 
This isomorphism gives an explicit presentation of higher 
K-groups modulo torsion via algebraic cycles. \\
Let us look more closely at the particular case of $K_1(X)$ for a smooth complex 
projective surface $X$. There it is known that
$CH^1(X,1)=\C^\times$ and $CH^p(X,1)=0$ for $p\ge 4$. The remaining interesting
parts of $K_1$ are therefore $CH^2(X,1)$ and $CH^3(X,1)$. The last group consists
of zero cycles on $X \times {\mathbb P}^1$ in good position and therefore
the map
$$
\tau: CH^2(X) \otimes_Z \C^\times \to CH^3(X,1), \quad x \otimes a
\mapsto (x,a) 
$$
is surjective. Therefore, the complexity of $CH^3(X,1)$ is governed by 
the complexity of $CH^2(X)$ which is fairly understood by Mumford's theorem
resp. Bloch's conjecture. We say that $CH^3(X,1)$ is decomposable. For $CH^2(X,1)$
the situation is quite different and the complex geometry of $X$ plays
an essential role in the understanding of it. The natural map
$$
\tau: CH^1(X) \otimes_Z \C^\times \to CH^2(X,1), \quad D \otimes a
\mapsto D \times \{a\} 
$$
is neither surjective nor injective in general. In the literature there
are several examples where the cokernel of $\tau$ is non-trivial modulo
torsion and even infinite dimensional, see \cite{C1}, \cite{GL},
\cite{SMS} and \cite{V1}. The kernel of $\tau$ is related but not equal to
$\Pic^0(X) \otimes \C^\times$ even modulo torsion by \cite[thm. 5.2]{S2}. Note that
the cokernel of $\tau$ is a {\sl birational invariant} 
(by localization) and hence 
vanishes for rational surfaces and, in fact, for all surfaces 
with geometric genus $p_g(X)=0$ and Kodaira dimension $\le 1$. The 
Bloch conjecture would imply that it vanishes also for all surfaces 
of general type which satisfy $p_g(X)=0$. One way to study 
$CH^2(X,1)$ is to look at the Chern class maps 
\eq{1.1}{
c_{2,1}: CH^2(X,1) \to H^3_\sD(X,\Z(2))=\frac{H^2(X,\C)}{H^2(X,\Z) +
F^2 H^2(X,\C)}.
}
The decomposable cycles (the image of $\tau$) are mapped to 
the subgroup
$$
{\rm NS}(X)\otimes_\Z \C^\times \subseteq 
\frac{H^2(X,\C)}{H^2(X,\Z) + F^2 H^2(X,\C)}
$$
generated by the {\sl N\'eron-Severi group} ${\rm NS}(X) \subset H^2(X,\Z)$
of all divisors in $X$. It is known \cite{SMS} that the image 
of $c_{2,1}$ is at most countable modulo this subgroup, so 
that the image of ${\rm coker}(\tau)$ in Deligne cohomology modulo
${\rm NS}(X)\otimes_\Z \C^\times $ is
at most countable. One conjectures that even ${\rm coker}(\tau)$ itself
is countable.\\
The Chern class  maps $c_{2,1}$ are defined as follows:
let $Z=\sum a_j Z_j \in CH^2(X,1)$ be a cycle. Each $Z_j$ is an integral curve
and inherits a rational map $f_j: Z_j \to \P^1$ from the projection map 
$X \times \P^1 \to \P^1$. Let $\gamma_0$ be a path on $\P^1$ connecting 
$0$ with $\infty$ along the real axis, then $\gamma:=\cup \gamma_j:= \cup
f_j^{-1}(\gamma_0)$ is a closed homological 1-cycle,
Poincar\'e dual to a cohomology
class in $F^2H^3(X,\Z)$ and therefore torsion, see \cite{SMS}. 
If we assume that $\gamma=0$ (for example if $b_1(X)=0$) then we write 
$\gamma=\partial \Gamma$ for a real piecewise smooth 2-chain $\Gamma$. 
The defining property of $c_{2,1}(Z)$ is, as a current, i.e. a linear functional on
differentiable complex valued 2-forms on $X$ 
\eq{1.2} { c_{2,1}(\alpha)= \frac{1}{2\pi i} 
\sum_j \int_{Z_j-\gamma_j} \log(f_j) \alpha + \int_\Gamma \alpha.
}
Now let $X$ be a projective K3 surface. Then $p_g(X)=1$, $b_2(X)=22$ 
and $b_1(X)=0$. The intersection form on $H^2(X,\Z)$ is known to be
the unimodular form $2E_8 \oplus 3 H$, where $H$ is the 2-dimensional
standard hyperbolic form. \\
The N\'eron-Severi lattice $NS(X) \subset H^{2}(X,\Z)$ 
has an orthogonal complement $T(X) \subset H^{2}(X,\Z)$. In particular
there is a well-defined morphism
$$
{\rm Coker}(\tau) \to \frac{T(X)\otimes \C^\times}{F^2}.  
$$
If we have an arbitrary smooth family $f:X \to B$ of complex algebraic 
surfaces over a quasiprojective complex variety $B$, and an algebraic family of 
cycles $Z_b \in CH^2(X_b,1)$ for all $b \in B$,
then we may define the {\sl normal function} 
$$
\nu(b):=c_{2,1}(Z_b) \in \frac{T(X_b)\otimes \C^\times}{F^2}.
$$
One can easily show that $\nu$ is a holomorphic (however multivalued) 
section of the corresponding
family of generalized tori $T(X_b)\otimes \C^\times/F^2$. 
Coming back to the case of K3-surfaces: there the canonical bundle $\omega_X$ 
is trivial, hence the group $H^{0,2}(X)=H^0(X,\Omega_X^2)^*=\C$ is 1-dimensional
and generated by the dual of $\omega_X$. In a smooth algebraic family $X_b$ of
K3-surfaces, the composition of the regulator with the projection onto 
$$
\frac{H^{0,2}(X_b)}{{\rm Im} H_2(X_b,\Z)}
$$
produces a multivalued holomorphic function on $B$, denoted by $\bar \nu(b)$,
which has poles at all $b$ where the 
family degenerates (proof see below). It is given by the formula
$$
\bar \nu(b)=\int_{\Gamma_b} \omega_{X_b},
$$
since the integral of $\omega_X$ over any effective divisor vanishes.
If $\sD_{\rm PF}$ denotes the Picard-Fuchs differential operator
of the Gau\ss-Manin connection associated to the family $X_b$ of K3-surfaces, then 
$\sD_{\rm PF}$ annihilates all periods of the family. Therefore we obtain the following 
result:
\begin{lem} Let $B \subset \bar B$ be smooth compactification of $B$. Then with the notation above, 
$\sD_{\rm PF}(\bar \nu)$ extends to a single-valued meromorphic function on $\bar B$ with poles only along 
degeneracies of $X_b$, and therefore satisfies a differential equation 
\eq{1.3}{
\sD_{\rm PF}(\bar \nu(b))=g(b),
}
where $g$ is a rational function in $b \in \bar B$. 
\end{lem}
The proof is given in the appendix.  Altogether we have obtained a map: 
\eq{1.4}{
\{ {\rm Families \; of \; Cycles \; in} \; CH^2(X_b,1) \}
\longrightarrow \{ {\rm Differential \; Equations/ Rational \; Functions} \}
}
For each such family of K3-surfaces it sends a family of cycles to the equation $\sD_{\rm PF}\bar\nu=g$ resp. 
the rational function $g$, which is the same information on a given family. One should view the
resulting solutions $\bar\nu(b)$ as {\it new transcendantal functions} arising from
the family of K-theoretic cycles in $CH^2(X_b,1)$. If $g$ is a non-trivial function, then $\bar \nu$ 
and hence $\nu$ is a non-flat section of the family of Deligne cohomology groups of $X_b$. In \cite{SMS}
the relationship between the infinitesimal behaviour of such normal 
functions and the mixed Hodge structure of the total space $X$ was already investigated. \\  
This situation is very reminiscent of a method developed
by Richard Fuchs \cite{F} in the case of the Legendre family $y^2=x(x-1)(x-t)$ of elliptic curves 
and investigated further in the work of Manin \cite[page 134]{M}. 
In particular there is a strong connection with differential 
equations of a generalized form of type Painlev\'e VI (loc. cit.). \\
There exists a {\it formula} to derive $g$: there is a so-called {\it inhomogenous Picard-Fuchs} equation 
\eq{1.5}{
\sD_{\rm PF} \omega_X = d_{rel} \beta
}
before integration over $\Gamma$, where $\beta$ is a section of the vector bundle
of (meromorphic) 1-forms in the fibers of the family $f:X \to B$. 
We say that $\Gamma$ {\it does not depend on} $b$ if it can be defined as real semi-algebraic subset
via {\it flat coordinates}, i.e. coordinate functions which are horizontal with respect
to the Gau\ss-Manin connection, and such that the defining inequalities of $\Gamma$ are polynomials 
not depending on $b$.  This shows on one hand that
for closed $\Gamma$ the periods satisfy the Picard-Fuchs equation, and on the other
hand for non-closed $\Gamma$ (not depending on $b$) with $\partial \Gamma = \gamma$ we get 
\eq{1.6}{
g(b)=\sD_{\rm PF} \int_\Gamma \omega_X = \int_\Gamma d_{rel} \beta = \int_\gamma \beta.
}
The last equality uses a version of Stokes theorem for currents since some of the
differential forms involved will in general have integrable singularities. Hence Stokes 
theorem for currents (see \cite[chap. 3]{GH}) also implies that $\beta$ is integrable over $\gamma$.   
In general $\Gamma $ depends also on $b$, and then there will be an additional
contribution from the derivatives of the boundaries of the integral.
In the case of the Legendre family $y^2=x(x-1)(x-t)$ of elliptic curves, 
$\beta$ is a meromorphic function ($0$-form)  
$$
\frac{y}{2(x-t)^2}, 
$$
by \cite[p. 310]{F},\cite[p. 76]{M}. Manin has 
put these equations into a more formal context (so-called $\mu$-equations) so that 
one can understand the sections and operators in a coordinate-free way in terms of certain
locally free sheaves on $B$. This plays also a role in his work on the functional Mordell conjecture.   
Furthermore, after uniformizing the elliptic curves, the
inhomogenous Picard-Fuchs equation is equivalent to a version involving the Weierstra\ss  
$\ \mathfrak{p}$-function
\cite[p. 137]{M}: 
$$
\frac{d^2 z}{d\tau^2} = \frac{1}{(2\pi i)^2} \sum_{j=0}^3 \alpha_j \mathfrak{p}_z(z+\frac{T_j}{2},\tau),
$$
where $\alpha_j$ are constants parametrizing the family of differential equations and 
$(T_0,\ldots,T_3)=(0,1,\tau,1+\tau)$ are the vertices of the fundamental parallelogram. In this way the 
transcendental aspect of the solutions and also the connection to integrable systems becomes apparent, 
see \cite[p. 139] {M}. In the future we hope to investigate further
the transcendental properties of our solutions (using again uniformization) and
study the attached integrable systems. \\
The rest of this article is devoted to a particular solution 
of the inhomogenous Picard-Fuchs equation for a certain family of K3 surfaces introduced in section 2.
In section 3 we deduce Esnault's conjecture from the non-vanishing of 
the $\sD_{\rm PF}(\bar \nu)$ in the special case $b=\sqrt{-1}$. In section 4 we study 
a certain Shioda-Inose model of $X_b$ which has isomorphic transcendental cohomology. 
This leads to an explicit computation of $\beta$ in this case. 

\section{An example: a mirror family of K3-surfaces}

We will study the one-dimensional family of K3-surfaces given by the quartic equations
\eq{2.1} { X_b:=\{(x,y,z,w) \in \P^3 \mid f(x,y,z,w)=xyz(x+y+z+bw)+w^4=0 \}.
}
with $b \in \P^1$. 
Note that this surface, for general $b$, 
is not smooth but has six singular points
defining a rational singularity of type $A_3$. The six points are
(\cite[sect.4]{NS}):
\begin{align*}
P_1= & (0,1,-1,0), & P_2= & (1,-1,0,0) \\
P_3= & (1,0,-1,0), & P_4 = & (1,0,0,0) \\
P_5= & (0,1,0,0), & P_6= & (0,0,1,0)
\end{align*}
The minimal resolution of the singularities defines a generically smooth family of K3-surfaces. 
In \cite{NS} the following theorem was shown:
\begin{thm}{\rm (Narumiya/Shiga)} The family $X_b$ has the following properties: 
\begin{enumerate} 
\item  It arises as a mirror family from the 
dual of the simplest polytope $P$ of dimension three. The dual mirror 
family is the family of all quartic K3-surfaces.
\item The rank of $\Pic(X)$ is $\ge 19$ for all $b \in \P^1 \setminus
\{0,\pm 4, \infty \}$ and equal to $19$ for very general $b$ (see loc.cit. \S 4.).
\item $T(X_b)$ has signature $(2,1)$ for $b \in \P^1 \setminus
\{0,\pm 4, \infty \}$. 
\item The periods of $X_b$ satisfy the Picard-Fuchs equation
$$(1-u)\Theta^3-\frac{3}{2}u\Theta^2-\frac{11}{16}u\Theta -
\frac{3}{32}u=0
$$
(where $\Theta=u\frac{d}{du}$)
of the generalized Thomae hypergeometric function \cite{T}
$$ F_{3,2}(\frac{1}{4}, \frac{2}{4},\frac{3}{4},1,1;u)
$$
and where we set $u:=(\frac{4}{b})^4$. 
\item In other words, the Picard-Fuchs equation is given by 
\eq{2.2}{ (1-u)u^2 \Phi''' + 3u(1-\frac{3}{2}u) \Phi'' 
+(1-\frac{51}{16}u)\Phi'  -\frac{3}{32}\Phi =0.
}
\item The mirror map of the family $X_b$ is given by the arithmetic 
Thompson series $T(q)$ of type {\bf 2A}in the classification of Conway and Norton \cite{CN}:
$$
T(q)=\frac{1}{q} +8+4372q + 96256 q^2+124002q^3 +10698752q^4 + \dots
$$
\end{enumerate}
\end{thm}
\proof We refer to \cite{NS} for more details, but we sketch the proof 
of (4) and (5) since this is crucial. (1)-(3) follow from the construction. 
In particular six $A_3$-singularities give rise to 18 independent cohomology
classes of type (1,1) so that the Picard number is $\ge 19$.
Since (5) is an easy corollary of (4), we prove (4). 
In \cite{NS} the periods are computed as power series in $1/b$ and the
differential equation in (4) follows from \cite{T}. 
As in \cite{NS} we consider the new affine equation 
$$
f(x,y,z)=xyz(x+y+z+1)+1/b^4=0
$$
obtained by substituting $w':=bw$ and setting $w'=1$. The periods are integrals 
of the form
$$
I(b)=\frac{1}{2\pi i} \int_{|x|=|y|=|z|=1/4} \frac{dxdydz}{xyz(x+y+z+1)+1/b^4}.
$$
On the other hand one has a geometric series expansion at $b=\infty$:
$$ 
\frac{1}{xyz(x+y+z+1)+1/b^4}= 
\sum_{n=0}^\infty \frac{(-1)^n b^{-4n}}{(xyz)^{n+1}(x+y+z+1)^{n+1}}.
$$
Changing the order of summation and integration, we obtain
\begin{align*} 
I(b)= & \frac{1}{2\pi i} \int_{|x|=|y|=|z|=1/4}\sum_{n=0}^\infty 
\frac{(-1)^n b^{-4n}dxdydz}{(xyz)^{n+1}(x+y+z+1)^{n+1}} \\
= & \frac{1}{2\pi i} \sum_{n=0}^\infty \int_{|x|=|y|=|z|=1/4} 
\frac{(-1)^n b^{-4n}dxdydz}{(xyz)^{n+1}(x+y+z+1)^{n+1}} 
\end{align*}
Now one can apply 3 times the residue theorem and gets
$$
I(b)=(2\pi i)^2 \sum_{n=0}^\infty \frac{(4n)!}{(n!)^4} b^{-4n}.
$$
Observing the identity involving Pochhammer symbols
$$
\frac{(4n)!}{(n!)^4} =
\frac{(\frac{1}{4})_n (\frac{2}{4})_n (\frac{3}{4})_n}{(1)_n (1)_n (1)_n} (4^4)^n,
$$
we have shown that
$$
I(b)=(2\pi i)^2 F_{3,2}(\frac{1}{4},\frac{2}{4},\frac{3}{4},1,1;(\frac{4}{b})^4)
$$
Substituting $u:=(\frac{4}{b})^4$, one gets a multiple of the functions
$F_{3,2}(\frac{1}{4},\frac{2}{4},\frac{3}{4},1,1;u)$ which satisfy a differential
equation of order three precisely of the type described in (4) resp. (5) 
by \cite{T}. \qed \\
\begin{cor} \label{pfcor} 
In $b$-coordinates, the Picard-Fuchs equation can be written as
\eq{2.3}{ ((\frac{b}{4})^4-1)(\frac{b}{4})^3 \Phi''' 
+ \frac{3}{4}(\frac{b}{4})^2 (1+(\frac{b}{4})^4) 
\Phi'' + \frac{1}{16} \frac{b}{4} ((\frac{b}{4})^4-6) \Phi' +\frac{3}{32}\Phi = 0  
}
\end{cor} 
\proof Use chain rule. \qed \\
\ \\
To make the following computations easier, we follow \cite{NS} and perform the following 
birational coordinate change (written in affine coordinates):
$$
X=xy, \quad Y=i(\frac{bxy}{2}+\frac{1+xy²z}{z}), \quad Z=yz.
$$
Then $X,Y,Z$ are affine coordinates and define the family of surfaces 
$S_b$ in Weierstra\ss form:
$$
S_b: \quad Y^2=X(X^2+X(Z+\frac{1}{Z}-\frac{b^2}{4})+1)
$$
as an elliptic fibration over $\P^1$ in the $Z$-coordinate. 
The inverse transformation is given by
$$
x=-2 \frac{iX(1+ZX)}{(-2Y+ibX)Z}, \quad y=\frac{1}{2}\frac{iZ(-2Y+ibX)}{1+ZX},
\quad z=-2 \frac{i(1+ZX)}{-2Y+ibX},
$$
in affine coordinates. The following is taken from \cite{NS} (with a slight correction):

\begin{lem} The surfaces $S_b$ are ramified coverings of $\P^1 \times \P^1$ 
(in $X,Z$ coordinates). In $ X,Y,Z$ coordinates, the canonical holomorphic two form on 
$S_b$ is given up to an non-zero constant by 
\eq{2.4}{
\omega=\frac{dXdZ}{YZ}=\frac{dXdY}{X^2(Z-\frac{1}{Z})},
}
where $X,Z$ are flat coordinates and 
$$
Y=Y(b)=\sqrt{P(X,Z)}=\sqrt{X(X^2+X(Z+\frac{1}{Z}-\frac{b^2}{4})+1)}.
$$
\end{lem}
\proof In the $x,y,z$ coordinate system the holomorphic two form is given up to a constant 
by $\frac{dxdy}{f_z}$ in affine $(x,y,z)$-coordinates with $w=1$. 
Then, using the coordinate transformations above, one computes that
$$
\omega=\sqrt{-1} \cdot \frac{dxdy}{f_z}= \frac{dXdY}{X^2(Z-\frac{1}{Z})}= \frac{dXdZ}{YZ},
$$
since $2YdY=X^2(1-\frac{1}{Z^2})dZ+ \frac{\partial P}{\partial X} dX$. \qed 
\ \\
Now, if we apply the Picard-Fuchs-Operator $\sD_{\rm PF}$ from Cor.\ref{pfcor} 
to $\frac{dXdZ}{YZ}$, we get an expression of the form 
$$
\sD_{\rm PF} \frac{dXdZ}{YZ}= K(X,Z) \cdot \frac{dXdZ}{Y^7 Z }
$$
where $K(X,Z)$ is a polynomial function in $X,Z$.

\section{The normal function and the Picard-Fuchs equation}

On any elliptic surface, the easiest way to find cycles in $K_1$ is to use fibers. 
However sometimes configurations coming from N\'eron fibers (degenerate into a union of 
$\P^1$'s) do have trivial class in $K_1$ as was already observed by Beilinson in \cite{B}.  
But one can use one smooth fiber together with a bunch of sections (rational) and rational 
curves in degenerate fibers. In our example let us take the following cycles: denote by
$S_b$ the surface defined by the equation 
\eq{3.1} {
S_b: \quad ZY^2=X(X^2Z+(Z^2+1-Z\frac{\; b^2}{4})X+Z)
}
Let $C_b$ be the smooth elliptic fibre over $Z=1$ of this surface. Its defining equation 
is hence
$$
C_b: \quad Y^2=X(X^2+(2-\frac{\; b^2}{4})X+1)
$$
The quadratic term $X^2+(2-\frac{b^2}{4})X+1$ in the right hand side
has two negative real roots if $b$ is purely imaginary, for example
$b=\sqrt{-1}$. 
The points $X=0$ and $X=\infty$ are rationally equivalent on $C_b$ 
after taking a multiple of two, since they are ramification points. 
The real line from $0$ to $\infty$ does not hit the
other ramification points by this observation. 
The surface $X_b$ in this birational model has $0$ and 
$X=\infty$ as sections. The fiber over $Z=0$ on $S_b$ decomposes into 
three rational curves with multiplicity counted. Hence one can construct a cycle in 
$CH^2(S_b,1)$
for general $b$ by using $C_b$, the two sections and the degenerate fibers and appropriate 
rational functions on all curves. In $X,Z$ coordinates the region $\Gamma$ is given by the
real square $0 \le Z \le 1, 0 \le X \le \infty$. For $b=\sqrt{-1}$ we make the following
observation: 
\begin{lem} For $b=\sqrt{-1}$ all coefficients occuring
in $\sD_{\rm PF} \frac{dXdZ}{YZ}= K(X,Z) \cdot \frac{dXdZ}{Y^7 Z }$ are positive integers, i.e. all 
coefficients of $K(X,Z)$ and all coefficients of $Y^2=X(X^2+X(Z+\frac{1}{Z}-\frac{b^2}{4})+1).$ 
\end{lem}
\proof Here the rules of differentiating are $\frac{\partial X}{\partial b}=
\frac{\partial Z}{\partial b}=0$ and $\frac{\partial Y}{\partial b}= \frac{b}{4}\frac{1}{Y^3}$. 
This implies that odd derivatives of $1/Y$ get multiplied by even powers of $b$. Now if we look at the
coefficients of equation \eqref{2.3}, we see that the coefficients at $\Psi'''$ and 
$\Psi'$ become positive, since  $(\frac{b}{4})^4-1$ and $(\frac{b}{4})^4-6$ are negative
rational numbers and get multiplied with $(\frac{b}{4})^6$, resp. $(\frac{b}{4})^2 $ which are
both also negative rational numbers. The coefficients at $\Psi$ and $\Psi''$ involve already 
4-th powers of $b$ and hence positive. Consequently all coefficients occuring are positive.   
Using any computer algebra program this can be verified and indeed one has:
\begin{align*}
\sD_{\rm PF} \frac{dXdZ}{YZ}= & {\displaystyle \frac {dXdZ}{8192}} (349951\,X^{3}\,
Z^{3} + 85952\,X\,Z^{3} + 85952\,X^{5}\,Z^{3} + 171904\,X^{4}\,Z
^{4} + 171904\,X^{4}\,Z^{2} \\
  & \mbox{} + 85952\,X^{3}\,Z^{5} + 171904\,X^{2}\,Z^{4} + 
294912\,X^{2}\,Z + 294912\,X\,Z^{2} + 98304\,Z^{3} \\
 & \mbox{} + 98304\,X^{3} + 294912\,X^{4}\,Z + 909352\,X^{2}\,Z
^{3} + 909352\,X^{3}\,Z^{2} + 98304\,X^{6}\,Z^{3} \\
  & \mbox{} + 294912\,X^{5}\,Z^{4} + 294912\,X^{5}\,Z^{2} + 
909352\,X^{4}\,Z^{3} + 294912\,X^{4}\,Z^{5} + 909352\,X^{3}\,Z^{4
} \\
  & \mbox{} + 294912\,X\,Z^{4} + 98304\,X^{3}\,Z^{6} + 294912\,X
^{2}\,Z^{5} + 85952\,X^{3}\,Z + 171904\,X^{2}\,Z^{2})/\\
  & ( 4\,X^{2}\,Z +4\,X\,Z^{2} +4\,X +X\,Z +4\,Z)^{7/2}\,
\sqrt{  {\displaystyle XZ} }
\end{align*}
This completes the proof. \qed \\
\ \\
In particular, if we integrate over the positive 
region $\Gamma$, we get a positive and non-zero integral. Since the boundary of $\Gamma$ 
is defined as the rectangle $0 \le Z \le 1, 0 \le X \le \infty$ and $X,Z$ are flat with 
respect to the connection, we say that $\Gamma$
does not depend on $b$ (see introduction) and this suffices to show that
the normal function is non-trivial. So we have proved the Esnault's conjecture (see \cite{SMS}):
\begin{cor} The projected normal function $\bar \nu(b)$ does not 
satisfy the Picard-Fuchs equation 
\eq{3.3}{ ((\frac{b}{4})^4-1)(\frac{b}{4})^3 \Phi''' 
+ \frac{3}{4}(\frac{b}{4})^2 (1+(\frac{b}{4})^4) 
\Phi'' + \frac{1}{16} \frac{b}{4} ((\frac{b}{4})^4-6) \Phi' +\frac{3}{32}\Phi = 0  
}
In particular, it is not a rational multiple of a period for
all but a countable number of values $b$.
For those $b$, the corresponding cycle $Z_b$ has no integer multiple 
which is decomposable modulo $Pic(X_b) \otimes {\mathbb C}^*$.  
\end{cor}
The main open problem remains to find 1-form $\beta$ such that $d \beta= \sD_{\rm PF} \frac{dXdZ}{YZ}$. 
We will compute such a $\beta$ for the Kummer type model of these K3 surfaces in the next section. 

\section{The solution of $\sD \omega=d\beta$}

In \cite{NS} you can find the description of a 2:1 map $\pi : S_b \to S_b'$ onto a Kummer
surface $S_b'$ which has a birational model with the equation
\eq{4.1}{ u^2=s(s-1)(s-(\frac{\nu+1}{\nu-1})^2)t(t-1)(t-\nu^2), 
}
where $\nu$ and $b$ are related via the algebraic equation
$$
b^2=-4 \cdot \frac{(\nu^2+1)^2}{\nu(\nu^2-1)}.
$$
We prefer to use this equation for a computation of the solution of $\sD \omega=d\beta$, since it is 
slightly easier but we do not loose essential information. 
In this description we see that the transcendental part of $H^2(S_b)$ also denoted by $T(S_b)$ 
has a Inose-Shioda structure (in the sense of \cite{Mo}) and is therefore related to a variation 
of a family of elliptic curves. In fact there are two  
isogenous elliptic curves $E_1(\nu)$ and $E_2(\nu)$ with equations
\eq{4.2}{
E_1(\nu): \; u_1^2=s(s-1)(s-\nu^2), \quad E_2(\nu): \; u_2^2=t(t-1)(t-(\frac{\nu+1}{\nu-1})^2)
}
together with a Nikulin involution (see \cite{VY}) on the abelian surface $A=E_1 \times E_2$
such that the associated Kummer surface is $S'_b$ and one has an isomorphism 
$T(S_b) \cong T(S_b')$ under $\pi_*$.  This explains in addition why the periods of $S_b$ 
are squares of other hypergeometric functions related to the family 
$E_1(\nu)$ resp. $E_2(\nu)$. More details about the birational map can be found in \cite{NS}. 
Further instances where Inose-Shioda structures and modular forms come up can be 
found in \cite{Do}. Let us now compute the Picard-Fuchs equation of the family $E_1(\nu)$. If we let $\lambda=\nu^2$, 
we have $\frac{\partial \nu}{\partial \lambda}=\frac{1}{2\nu}=\frac{1}{2\sqrt{\lambda}}$ 
and therefore for any function $\Phi$ we have the transformation rules
$$
\frac{\partial \Phi}{\partial \lambda}= 
\frac{\partial \Phi}{\partial \nu} \cdot  \frac{\partial \nu}{\partial \lambda}=
\frac{\partial \Phi}{\partial \nu} \frac{1}{2\sqrt{\lambda}}=
\frac{\partial \Phi}{\partial \nu} \frac{1}{2\nu}
$$
and for the second derivative
$$ 
\frac{\partial^2 \Phi}{\partial \lambda^2}=\frac{1}{4\nu^2} \frac{\partial^2 \Phi}{\partial \nu^2}-
\frac{1}{4\nu^3} \frac{\partial \Phi}{\partial \nu}.
$$
Plugging this into the standard hypergeometric Picard-Fuchs equation
$$
\lambda(1-\lambda)\Phi''(\lambda)+(1-2\lambda)\Phi'(\lambda)-\frac{1}{4} \Phi(\lambda)=0,
$$
we get the new equation
$$
(1-\nu^2)\Phi''(\nu)+\frac{1-3\nu^2}{\nu} \Phi'(\nu) -\Phi(\nu) =0
$$
and the inhomogenous variant (equality of 1-forms)
$$ 
(1-\nu^2)  \frac{\partial^2}{\partial \nu^2} \omega(s)   
+\frac{1-3\nu^2}{\nu}  \frac{\partial}{\partial \nu} \omega(s)  - \omega(s) 
= 2 d_{\rm rel} \frac{\sqrt{s(s-1)(s-\nu^2)}}{(s-\nu^2)^2}
=  2 d_{\rm rel} \frac{s^2(s-1)^2}{\sqrt{s(s-1)(s-\nu^2)}^3},
$$
where
$$ 
\omega(s)= \frac{ds}{\sqrt{s(s-1)(s-\nu^2)}} .
$$
The rational normalized version of this equation is
\eq{4.3}{ 
\frac{\partial^2}{\partial \nu^2} \omega(s)   
+\frac{1-3\nu^2}{\nu(1-\nu^2)}  \frac{\partial}{\partial \nu} \omega(s)  - \frac{1}{1-\nu^2} \omega(s) 
= \frac{2}{1-\nu^2} d_{\rm rel} \frac{\sqrt{s(s-1)(s-\nu^2)}}{(s-\nu^2)^2},
}
In a similar way we use the substitution $\lambda=(\frac{\nu+1}{\nu-1})^2$ and get the 
formula $\frac{\partial \nu}{\partial \lambda}= -\frac{(\nu-1)^3}{4(\nu+1)}$ and hence 
$$
\frac{\partial \Phi}{\partial \lambda}=
\frac{\partial \Phi}{\partial \nu} \cdot  \frac{\partial \nu}{\partial \lambda}=
 -\frac{(\nu-1)^3}{4(\nu+1)}\frac{\partial \Phi}{\partial \nu}, \quad 
\frac{\partial^2 \Phi}{\partial \lambda^2}= 
\frac{(\nu-1)^6}{16(\nu+1)^2} \frac{\partial^2 \Phi}{\partial \nu^2} + 
\frac{(\nu-1)^5(\nu+2)}{8(\nu+1)^3} \frac{\partial \Phi}{\partial \nu}
$$
Combining all this, we get the equation 
\eq{4.4} {
\frac{\partial^2}{\partial \nu^2} \omega(t)   
+\frac{\nu^2-2\nu-1}{\nu(\nu^2-1)}  
\frac{\partial}{\partial \nu} \omega(t)  + \frac{1}{\nu(\nu-1)^2} \omega(t) 
= -\frac{2}{\nu(\nu-1)^2} d_{\rm rel} \frac{\sqrt{t(t-1)(t-(\frac{\nu+1}{\nu-1})^2}}
{(t-(\frac{\nu+1}{\nu-1})^2)^2}
} 
for 
$$
\omega(t)= \frac{dt}{\sqrt{t(t-1)(t-(\frac{\nu+1}{\nu-1})^2 )}} .
$$
We have to compute a sort of convolution product of these two equations in the following sense:
set 
$$
\omega=\omega(s) \wedge \omega(t)=
\frac{dsdt}{\sqrt{s(s-1)(s-\nu^2)t(t-1)(t-(\frac{\nu+1}{\nu-1})^2})},
$$
and notice that we have the product formula:
$$
\frac{\partial^3}{\partial \nu^3} \omega = \frac{\partial^3}{\partial \nu^3} \omega(s) 
\wedge \omega(t)+ 
3 \frac{\partial^2}{\partial \nu^2} \omega(s) \wedge \frac{\partial}{\partial \nu} \omega(t) +
 3 \frac{\partial}{\partial \nu} \omega(s) \wedge \frac{\partial^2}{\partial \nu^2} \omega(t)+
\omega(s) \wedge  \frac{\partial^3}{\partial \nu^3} \omega(t),
$$
Similar formulas hold for lower derivatives. Note that $s,t$ are flat coordinates so that
differentiating a differential form with respect to the coefficients is a well-defined procedure.
Such formulas can be used to compute
$\frac{\partial^3}{\partial \nu^3} \omega$ and obtaining a Picard-Fuchs differential operator $\sD$ 
for $\omega$ together with a solution $\beta$ of $\sD \omega=d_{\rm rel}\beta$:  
\begin{lem} One has the following inhomogenous Picard-Fuchs equation involving 2-forms:
\eq{4.5}{ 
\frac{\partial^3\omega}{\partial \nu^3}  + 3 \frac{2\nu+1}{\nu(\nu+1)}
\frac{\partial^2\omega}{\partial \nu^2}  + 
\frac{7\nu^4-6\nu^3-4\nu^2+6\nu+1}{(\nu-1)^2(\nu+1)^2\nu^2}
\frac{\partial\omega}{\partial \nu} 
+ \frac{\nu^4-2\nu^3-2\nu-1}{(\nu-1)^3(\nu+1)^2\nu^2} \omega = d_{\rm rel} \beta,
}
where
\eq{4.6} {
\beta = -2 \frac{s(s-1)(2\nu^4+3\nu^3-\nu^2-3\nu s+s \nu^2-2s)}{\nu(s-\nu^2)^2(\nu^2-1)^2 
\sqrt{s(s-1)(s-\nu^2)}}
\omega(t)+ \frac{6}{1-v^2} \frac{\sqrt{s(s-1)(s-\nu^2)}}{(s-\nu^2)^2} \omega'(t) 
}
$$
2 \frac{t(t-1)(2t\nu^4-7t\nu^3+7t\nu^2-t\nu-t-2\nu^4-7\nu^3-7\nu^2-\nu+1)}
{\nu^2(\nu-1)^4(t-(\frac{\nu+1}{\nu-1})^2)^2(\nu^2-1) 
\sqrt{t(t-1)(t-(\frac{\nu+1}{\nu-1})^2})} \omega(s)  
+6 \frac{\sqrt{t(t-1)((t-(\frac{\nu+1}{\nu-1})^2)}}{\nu(\nu-1)^2(t-(\frac{\nu+1}{\nu-1})^2)^2} 
\omega'(s).
$$
\end{lem}
\proof We denote the derivative of a function (or a form in flat coordinates) $f$ of $\nu$ 
by $f'$. Let us carry out the computation in a more general setting: assume
we have two Picard-Fuchs equations:
$$
\omega''(s) -A_s \omega'(s) - B_s \omega(s) = d\beta_s,
$$
$$
\omega''(t) - A_t \omega'(t) - B_t \omega(t) = d\beta_t,
$$ 
with functions $\beta_s,\beta_t$ and $A_s,B_s,A_t,B_t$ depending on $\nu$.  Now first note that
$$
\omega'''(s)=A_s \omega''(s)+(A_s'+B_s) \omega'(s) +B_s' \omega(s) + \frac{d}{d\nu} d\beta_s
$$
and furthermore
$$
\frac{d}{d\nu} d\beta_s=d \beta_s', 
$$
by symmetry of mixed derivatives. A similar relation holds for $t$. By using the product formulas
$$
\omega''' = \omega'''(s) \wedge \omega(t)+ 3 \omega''(s) \wedge \omega'(t) +
 3 \omega'(s) \wedge \omega''(t)+ \omega(s) \wedge  \omega'''(t),
$$
$$
\omega'' = \omega''(s) \wedge \omega(t)+ 2 \omega'(s) \wedge \omega'(t) 
+ \omega(s) \wedge \omega''(t),
$$
we obtain that
$$
\omega'''-\frac{3}{2}(A_s+A_t)\omega''
=[A_s'+B_s+3B_t-\frac{1}{2} A_s^2 - \frac{3}{2} A_s A_t] 
\omega'(s) \wedge \omega(t)
$$
$$
+[A_t'+B_t+3B_s-\frac{1}{2} A_t^2 - \frac{3}{2} A_s A_t] \omega(s) \wedge \omega'(t)
+[A_sB_s+A_tB_t-\frac{3}{2}(B_s+B_t)(A_s+A_t)+B_s'+B_t'] \omega
$$
$$
+d_{\rm rel}[(\beta_s' -(\frac{1}{2} A_s + \frac{3}{2} A_t)\beta_s) \omega(t) 
+3 \beta_s \omega'(t) 
- (\beta_t'-(\frac{1}{2} A_t + \frac{3}{2} A_s)\beta_t) \omega(s) -3 \beta_t \omega'(s)]
$$ 
does not involve anymore terms of the form $\omega'(s) \wedge \omega'(t)$. Now let
$$
A:=\frac{3}{2}(A_s+A_t), \; B:= -\frac{1}{2}A_s^2 -\frac{3}{2}A_sA_t+A_s'+B_s+3B_t,
$$
$$ 
C:=-\frac{1}{2} A_sB_s - \frac{3}{2}(A_sB_t+A_tB_s) -\frac{1}{2} A_tB_t +B_s' +B_t'.
$$
Then, assuming that we have the equality
$$
-\frac{1}{2}A_s^2-\frac{3}{2}A_sA_t +A_s'+B_s+3B_t=
-\frac{1}{2}A_t^2-\frac{3}{2}A_sA_t +A_t'+B_t+3B_s,
$$
(this condition is equivalent to the fact that the elliptic curves
$E_1(\nu)$ and $E_2(\nu)$ are isogenous), then we have the following  
inhomogenous Picard-Fuchs equation: 
$$
\omega''' - A \omega'' -B \omega' - C \omega = d_{\rm rel} \beta,
$$
where $\beta$ is the 1-form: 
$$
\beta:= (\beta_s' -(\frac{1}{2} A_s + \frac{3}{2} A_t)\beta_s)) \omega(t) +3 \beta_s \omega'(t) 
+ (\beta_t'-(\frac{1}{2} A_t + \frac{3}{2} A_s)\beta_t)) \omega(s) +3 \beta_t \omega'(s). 
$$ 
In our case $A_s=-\frac{1-3\nu^2}{\nu(1-\nu^2)}$, $B_s=\frac{1}{1-\nu^2}$,
$A_t=-\frac{\nu^2-2\nu-1}{\nu(\nu^2-1)}$ and $B_t=-\frac{1}{\nu(\nu-1)^2}$. Therefore we get
$$
A= -3 \frac{2\nu+1}{\nu(\nu+1)}, 
\quad B=-\frac{7\nu^4-6\nu^3-4\nu^2+6\nu+1}{(\nu-1)^2(\nu+1)^2\nu^2}, \quad
C= - \frac{\nu^4-2\nu^3-2\nu-1}{(\nu-1)^3(\nu+1)^2\nu^2}
$$
and for $\beta$ the expression above. This finishes the proof. \qed

\section{Appendix}

In this section we give the proof of the following lemma from the introduction: \\
\ \\
{\bf Lemma 1.1:} {\sl In the situation of the introduction, $\sD_{\rm PF}(\bar \nu)$
is a single-valued meromorphic function on $B$ with poles only along 
degeneracies of $X_b$, and therefore satisfies a differential equation 
$$
\sD_{\rm PF}(\bar \nu)=g,
$$
where $g$ is a rational function in $b \in B$. } \\
\noindent \proof Assume that we have a family $f:\bar X \to \bar B$ 
of projective surfaces over a compact Riemann surface $\bar B$. 
Let $\Sigma \subseteq \bar B$ be the finite
subset over which there are singular fibers. Let $h: X \to B$ be the smooth part of $f$. 
We may assume that the family is semi-stable (semi-stable reduction) 
and that there is a cycle $Z \in CH^2(X,1)$ such that the restriction
of $Z$ to all fibers induces the family of cycles in $CH^2(X_b,1)$ (both reductions 
require perhaps a finite cover of $B$ which does not however change the assertion).    
The cycle $Z$ has a class $c_{3,2}(Z) \in H^3_\sD(X,\Z(2))$ in Deligne cohomology.
By semi-stability $\Delta:=f^{-1}\Sigma$ is a divisor with
strict normal crossings and its Deligne cohomology can be computed via the logarithmic 
de Rham complex. Let $\sV^2$ be the sheaf of transcendental cohomology classes in 
$R^2h_*\C$, a local system of rank $22-\rho(X_b)$ for $b$ general. The Deligne class vanishes in 
$F^3 \cap H^3(X_b,\Z)$, since $b_3(X_b)=0$ in our case, and therefore
induces a holomorphic normal function $\nu \in H^0(B,\sV^2 \otimes \sO_B/F^2)$ over $B$.
However since the family is semi-stable, there is a canonical extension of $\nu$ to a 
holomorphic section of the sheaf 
$R^2 f_* \Omega^*_{\bar X/\bar B}(log \Delta)$. 
This can be seen as follows: let 
$$
\Z_{\sD,X}(2)=Cone(Rj_*\Z(2)  \to 
\Omega^*_{\bar X}(log \Delta)/F^2)[-1]
$$ 
be the Beilinson-Deligne complex \cite{EV} of 
$X$ using the inclusion $X {\buildrel j \over \hookrightarrow} \bar X$ and 
$$
\Z_{\sD,h}(2)=Cone(Rj_*\Z(2) \to 
\Omega^*_{\bar X/\bar B}(log \Delta)/F^2)[-1]
$$ 
the relative Beilinson-Deligne complex of $h$. There is a natural surjection of complexes 
$\Z_{\sD,X}(2)\to \Z_{\sD,h}(2)$ which induces a morphism 
$H^3_\sD(X,\Z(2)) \to H^0(\bar B,R^3 f_*\Z_{\sD,f}(2))$. Since all fibers of $h$
satisfy $b_3(X_b)=0$, we conclude that the image of this element in
$H^0(\bar B,R^3 f_* Rj_*\Z(2))$ vanishes. Therefore 
the image of $c_{3,2}(Z)$ in $H^0(\bar B,R^3 f_*\Z_{\sD,f}(2))$ is coming (at least locally
because of monodromy) from a class in 
$H^0(\bar B, R^2 f_* \Omega^*_{\bar X/\bar B}(log \Delta)/F^2)$
and is thus an extension of $\nu \in H^0(B, R^2 f_* \Omega^*_{X/B}/F^2)$ to $\bar B$. 
By construction it is meromorphic along $\Sigma$ but still multivalued with indeterminacies 
in the local system of integral cohomology classes. Now we apply the 
Picard-Fuchs operator. This makes $g(b)$ a single-valued complex function on $B$.  
$\sD_{\rm PF}$  has meromorphic (rational) coefficients in $b$, since they are the coefficients of
the characteristic polynomial of the Gau\ss-Manin connection, which has regular singular points along 
$\Sigma$ by Deligne \cite{Del}. Therefore the resulting function $g(b)$ is holomorphic outside $\Sigma$, but can have 
higher order poles along $\Sigma$. By Chow's theorem any meromorphic function on $\bar B$ is rational. \qed  \\
\ \\ 
{\bf Acknowledgement:} We are grateful to S. Bloch, H. Esnault, Y. Manin and
J. Nagel for several valuable discussions in which they urged us to find the
solution of the inhomogenous equation. We are also grateful to N. Narumiya, H. Shiga
and N. Yui for explaining their results and to M. Saito for explaining his conjectures
about the structure of $CH^2(X,1)$ and his independent work on the same problem. We are also 
grateful to the referee for his suggestions to improve the paper.

\bibliographystyle{plain}
\renewcommand\refname{References}

\end{document}